\documentclass[12pt]{article}  
\usepackage{hyperref}
\usepackage{amssymb}
\let\oldmarginpar\marginpar
\renewcommand\marginpar[1]{${}^\clubsuit$\oldmarginpar[\raggedleft\scriptsize\sf #1]{\raggedright\scriptsize\sf #1}}
\newcommand{\cref}[1]{Corollary~\ref{#1}}
\newcommand{\dref}[1]{Definition~\ref{#1}}

\newcommand{\eqref}[1]{(\ref{#1})}
\newcommand{\lref}[1]{Lemma~\ref{#1}}
\newcommand{\pref}[1]{Proposition~\ref{#1}}

\newcommand{\tref}[1]{Theorem~\ref{#1}}
\newcommand{\Ee}{\mathbb {E}}
\newtheorem{theorem}{Theorem}

\newtheorem{proposition}[theorem]{Proposition}
\newtheorem{lemma}[theorem]{Lemma}
\newtheorem{corollary}[theorem]{Corollary}
\newtheorem{remark}[theorem]{Remark}

\newtheorem{definition}[theorem]{Definition}

\newcommand{\R}{\mathbb{R}}

\newcommand{\Le}{\mathbb{L}}
\newcommand{\Ve}{\mathbb{V}}

\newcommand{\spa}{\mbox{span}}
\newcommand{\hess}{\mbox{Hess\,}}

\newcommand{\Ima}{\mbox{Im }}
\newcommand{\rank}{\mbox{rank\ }}
\newcommand{\kerl}{\mbox{Ker }}

\newcommand{\po}{{\hspace*{-1ex}}{\bf .  }}
\newcommand{\ii}{isometric immersion }
\newcommand{\iis}{isometric immersions }

\newcommand{\N}{{\cal N}}
\newcommand{\Tes}{{\cal T}}
\newcommand\fall{\;\;\;\;\mbox{for all}\;\;}

\def\Sal{{\cal S}}

\def\Les{{\cal L}}

\def\Ies{{\cal J}}

\def\Kes{{\cal K}}

\def\<{\langle}

\def\>{\rangle}
\def\a{\alpha}
\def\va{\varphi}

\def\bea{\begin{eqnarray*} }
\def\eea{\end{eqnarray*} }
\def\be{\begin{equation} }
\def\ee{\end{equation} }

\def\nap{\nabla^\perp}

\def\nab{\widetilde\nabla }
\def\proof{\noindent{\it Proof: }}
\def\qed{\ifhmode\unskip\nobreak\fi\ifmmode\ifinner
\else\hskip5 pt \fi\fi\hbox{\hskip5 pt \vrule width4 pt
height6 pt  depth1.5 pt
\hskip 1pt }}
\begin{document}

\title{Genuine deformations of submanifolds II:\\the conformal case}  
\author {Luis A. Florit \& Ruy Tojeiro}
\date{}
\maketitle
\addtocounter{equation}{1}

\begin{abstract}  
We extend to the conformal realm the concept of genuine deformations of
submanifolds, introduced by Dajczer and the first author for the
isometric case. Analogously to that case, we call a conformal
deformation of a submanifold $M^n$ genuine if no open subset of $M^n$ can
be included as a submanifold of a higher dimensional conformally
deformable submanifold in such a way that the conformal deformation of
the former is induced by a conformal deformation of the latter. We
describe the geometric structure of a submanifold that admits a genuine
conformal deformation and give several applications showing the unifying
character of this concept.
\end{abstract}

\medskip
\section{Introduction}  
An \ii $\hat f\colon\,M^n\to\R^{n+q}$ with codimension $q$ of an
$n$-dimensional Riemannian manifold $M^n$ into Euclidean space is
said to be a {\it genuine isometric deformation\/} of a given \ii
$f\colon\, M^n\to\R^{n+p}$ if $f$ and $\hat f$ are nowhere (i.e.,
on no open subset of $M^n$) compositions, $f=F\circ j$ and $\hat
f=\hat F\circ j$, of an isometric embedding \mbox{$j\colon\,
M^n\hookrightarrow N^{n+r}$}  into a Riemannian manifold $N^{n+r}$
with $r>0$ and \iis $F\colon\,N^{n+r}\to\R^{n+p}$ and $\hat
F\colon\,N^{n+r}\to\R^{n+q}$:

\bigskip

\begin{picture}(110,84)
\put(155,31){$M^n$}
\put(209,31){$N^{n+r}$}
\put(413,31){$(1)$}
\put(242,62){$\R^{n+p}$}
\put(242,0){$\R^{n+q}$}
\put(195,59){$f$}
\put(195,4){${\hat f}$}
\put(240,46){${}_F$}
\put(239,26){${}_{\hat F}$}
\put(193,40.5){${}_j$}
\put(228,42){\vector(1,1){16}}
\put(228,28){\vector(1,-1){16}}
\put(177,44){\vector(3,1){60}}
\put(177,26){\vector(3,-1){60}}
\put(185,34){\vector(1,0){21}}
\put(185,36){\oval(7,4)[l]}
\end{picture}
\bigskip\medskip

\noindent More geometrically, an isometric deformation of an Euclidean
submanifold $M^n$ is genuine if no open subset of $M^n$ can be included
as a submanifold of a higher dimensional isometrically deformable
submanifold in such a way that the isometric deformation of the former is
induced by an isometric deformation of the latter.

This concept was introduced in \cite{df1}, where it was proved that if an
\ii $f\colon\, M^n\to\R^{n+p}$ and a genuine isometric deformation
$\hat f\colon\,M^n\to\R^{n+q}$ of it have sufficiently low codimensions
then they are {\it mutually (isometrically) ruled}, that is, $M^n$
carries an integrable \mbox{$d$-dimensional} distribution $D^d\subset TM$
whose leaves are mapped diffeomorphically by $f$ and $\hat f$ onto open
subsets of affine subspaces of $\R^{n+p}$ and $\R^{n+q}$,
respectively. The authors also obtained a sharp estimate on the dimension
$d$ of the rulings and proved that the normal connections and second
fundamental forms of $f$ and $\hat f$ satisfy strong additional
relations.

Besides containing several previous results on isometric deformations of
submanifolds as particular cases, this concept has given new geometric
insight on the structure of isometrically deformable submanifolds,
showing that genuinely deformable submanifolds are rather special and
providing an important step for extending to higher codimensions the
classical Sbrana-Cartan theory of isometrically deformable hypersurfaces
(\cite{sb}, \cite{ca1}, \cite{dft}).

\vskip .3cm

Our goal in this article is twofold. First, to extend the notion of
genuine deformations to the conformal realm, and to give a similar
description as in \cite{df1} of the geometric nature of a submanifold
that admits such a deformation. In particular, to provide a unified
account of several known results on conformal deformations of
submanifolds. Second, to understand geometrically the similitude between
the theories of isometric and conformal deformations of submanifolds. In
order to state our results we first set up some terminology.

A {\em conformal structure\/} on a manifold $M^n$ is an equivalence class
of conformal Riemannian metrics on $M^n$. Recall that two
Riemannian metrics $\<\,,\,\>$ and $\<\,,\,\>'$ on $M^n$ are
conformal if there exists a smooth function $\va$ on $M^n$ such that
$\<\,,\,\>'=\va^2\<\,,\,\>$. We call $\va$ the {\it conformal factor}
relating the metrics $\<\,,\,\>$ and $\<\,,\,\>'$.
Clearly, every Riemannian manifold has a canonical conformal structure
determined by its metric.

Given an immersion $f\colon\,M^n\to \bar M^m$ between differentiable
manifolds, since conformal metrics on $\bar M^m$ are pulled-back by $f$ to
conformal metrics on $M^n$, a conformal structure on $\bar M^m$ induces a
conformal structure on $M^n$, the {\em conformal structure on $M^n$
induced by $f$\/}. If $M^n$ is already endowed with a conformal
structure, we call $f$ {\em conformal\/} if such conformal structure
coincides with that induced by $f$.

A pair $\{f,\bar f\}$ of conformal immersions $f\colon\, M^n\to \R^{n+p}$
and $\bar f\colon\,M^n\to \R^{n+q}$ will be referred to simply
as a {\em conformal pair\/}. We say that the conformal pair
$\{f,\bar f\}$ {\it extends conformally} when there exist a conformal
embedding $j: M^n \to N^{n+r}$, with $r\geq 1$, and conformal immersions
$F:N^{n+r} \to \R^{n+p}$ and $\bar F:N^{n+r} \to \R^{n+q}$ such that
$f=F\circ j$ and $\bar f=\bar F\circ j$; see (1). We call the (ordered)
conformal pair $\{F,\bar F\}$ a {\it conformal extension} of
$\{f,\bar f\}$.

The conformal pair $\{f,\bar f\}$ is said to be {\it genuine} if there is
no open subset $U \subseteq M^n$ such that the restricted pair
$\{f|_U,\bar f|_U\}$ extends conformally. If $\{f,\bar f\}$ is a genuine
conformal pair, we also say that each of its elements is a {\em genuine
conformal deformation\/} of the other.

A conformal immersion $f\colon\,M^n\to\R^{n+p}$ is {\it genuinely
conformally rigid} in $\R^{n+q}$ for a fixed integer $q> 0$
if, for any given conformal immersion $\bar{f}\colon\,M^n\to \R^{n+q}$,
there is an open dense subset $U\subset M^n$ such that the pair
$\{f|_U,\bar{f}|_U\}$ extends conformally.

Similar definitions can be given for any ambient spaces carrying
conformal structures, as well as for {\em isometric\/} immersions between
arbitrary {\em semi\/}-Riemannian manifolds, in the same way as in
\cite{df1}.

\vspace{.2cm}

We say that an immersion $f\colon\,M^n\to\R^{n+p}$ is
$D^d$-{\em conformally ruled\/} if $M^n$ carries an integrable
\mbox{$d$-dimensional} distribution $D^d\subset TM$ whose leaves are
mapped diffeomorphically by $f$ onto open subsets of affine subspaces
or round spheres of $\R^{n+p}$. Then, at each point $x\in M$ we have a
symmetric bilinear form
$\beta^f=\beta^f(x)\colon\,T_xM\times T_xM\to T_x^\perp M$ defined by
$$
\beta^f(Z,X):=\a^f(Z,X)-\<\,Z,X\>\eta(x)
$$
and a subspace of the normal space $T_x^\perp M$ of $f$ at $x$ given by
$$
L^c_D(x)=L^c_D(f)(x)=\spa\{\beta^f(Z,X) : Z\in D^d(x)\;
\,\mbox{and}\,\; X\in T_xM\}.
$$
Here $\a^f\colon\,TM\times TM\to T_f^\perp M$ stands for the second
fundamental form of $f$ and $\eta(x)$ for the normal component of the
mean curvature vector of the (image by $f$ of the) leaf of $D$ through
$x\in M^n$. We always work on open subsets where the dimension of
$L^c_D(x)$ is constant, in which case such subspaces form a smooth
subbundle of $T_f^\perp M$ that we denote by $L^c_D$.  Observe that
$D\subseteq \N(\beta^f_{{L^c_D\!\!}^\perp})$, where $\N(\beta)$ denotes
the nullity space of a symmetric bilinear form $\beta$, and a subspace
as a subscript means to take the orthogonal projection onto that subspace.

We are now in a position to state the main result of this paper. As we
will see, it implies or even generalizes main results in several other
works, e.g. \cite{cd}, \cite{df4}, \cite{dt1}, \cite{dt2}, \cite{dt3},
whereas revealing their true geometric nature.

\begin{theorem}\label{int}\po
Let $f\colon M^n\!\to\R^{n+p}$ and $\bar{f}\colon M^n\to\R^{n+q}$ form a
genuine conformal pair, with $p+q\leq n-3$ and $\min\,\{p,q\}\le 5$.
Then, along each connected component of an open dense subset of $M^n$,
the immersions $f$ and $\bar{f}$ are mutually conformally $D^d$-ruled,
with
$$
d\geq n-p-q+3\ell^c_D,
$$
and $D^d=\N(\beta^{f}_{L^\perp})\cap \N(\beta^{\bar f}_{\bar
L^\perp})$, where $L:=L^c_D(f)$, $\bar L:=L^c_D(\bar f)$ and
$\ell^c_D:=\rank L=\rank \bar L$. Moreover, there exists a
parallel vector bundle isometry $\Tes:L\to\bar L$ such that
$\beta^{\bar f}_{\bar L}=\va\Tes\circ \beta^f_L$, where $\va$ is
the conformal factor relating the metrics induced by $f$ and $\bar
f$.
\end{theorem}

In other words,
up to an identification, the normal bundles of the immersions contain a
subbundle $L=\bar L$ with the same normal connections and the same
(conformal) second fundamental forms. On the other hand, the common
conformal rulings $D^d$ of the immersions are the  nullity of the
(conformal) second fundamental forms on their orthogonal complements
$L^\perp$ and $\bar L^\perp$. The larger is $L$, the bigger is $d$.
We point out that Example~2 in \cite{df1} also shows that the
estimate on $d$ in \tref{int} is sharp.

\vskip .2cm

As a consequence of Theorem \ref{int}, we obtain the following conformal
version of the main result of \cite{df3}. We recall from \cite{cd} that
the {\em conformal $s$-nullity\/} $\nu^c_s(x)$ of an immersion
$f\colon\,M^n\to\R^{n+p}$ at a point $x\in M^n$ is defined for
$1\leq s\leq p$ by
$$
\nu^c_s(x)= \max \{\dim\,\N\left(\a^f_{V^s}(x)-\<\,,\,\>_f\zeta\right)
\colon\,V^s\subset T_{x}^\perp M,\,\,\zeta\in V_s\},
$$
where $\<\,,\,\>_f$ stands for the metric on $M^n$ induced by $f$.

\begin{corollary}\po\label{cor1}
Let $f\colon\,M^n\to\R^{n+p}$ be an immersion and let $q$ be a positive
integer with $p\leq q\leq n-p-3$. Suppose that $p\leq 5$ and that $f$
satisfies
$$
\nu^c_s\leq n+p-q-2s-1\ \fall\ \ 1\leq s\leq p.
$$
For $q\geq p+5$ assume further that $\nu^c_1\leq n-2(q-p)+1$. Then, any
immersion \mbox{$\bar f\colon\,M^n\to\R^{n+q}$} conformal to $f$ is locally a
composition, i.e., there exists an open dense subset $V\subseteq M^n$
such that the restriction $\bar f$ to any connected component
$U$ of $V$ satisfies $\bar f|_U=h\circ f|_U$, where
$h\colon\,W\subset\R^{n+p}\to\R^{n+q}$ is a conformal immersion of an
open subset $W\supset f(U)$.
\end{corollary}

For $p=q$, the preceding corollary extends up to codimension $p=5$ the
main theorem of \cite{cd}, which ensures conformal rigidity of $f$ in
$\R^{n+p}$ whenever $p\leq 4$ and $\nu^c_s\leq n-2s-1$ for all
$1\leq s\leq p.$ The latter, in its turn, is a generalization of Cartan's
classical criterion $\nu^c_s\leq n-2$ for conformal rigidity of
hypersurfaces. Corollary
\ref{cor1} also generalizes the main result of \cite{dt2}, which deals
with the special case $p=1$, as well as \cite{s} up to codimension~$5$.

If we apply Theorem~\ref{int} for $p=q=2$, the estimate on $d$
implies that $\ell^c_D\le 1$.  This yields $d\geq n-4$ if $\ell^c_D=0$
and $d\geq n-1$ if $\ell^c_D=1$. In both cases, the conformal
nullity $\nu^c=\nu^c_2$ of both immersions satisfies $\nu^c\geq n-4$.
Therefore, under the assumptions that $n\geq 7$ and $\nu^c_f\le n-5$
everywhere, we conclude that $f$ is genuinely conformally rigid. Here,
this means that if $\bar f$ is a conformal deformation of $f$, then there
exists an open and dense subset of $M^n$ on each connected component of
which either $\bar f$ is conformally congruent to $f$ or $f$ can be
included as a submanifold of either a conformally flat or a Cartan
hypersurface and $\bar f$ is induced by a conformal deformation of such
hypersurface. This is the content of the main theorem in \cite{dt3}.

\vskip .2cm

Clearly, Theorem \ref{int} yields the following criterion for genuine
conformal rigidity.

\begin{corollary}\po\label{rig1}
Let $f\colon\, M^n\to\R^{n+p}$ be a conformal immersion
and let $q$ be a positive integer with $p+q\leq n-3$ and
$\min\,\{p,q\}\le 5$. If $f$ is not $(n\!-\!p\!-\!q)$--conformally ruled
on any open subset of $M^n$, then $f$ is genuinely conformally rigid in
$\R^{n+q}$.
\end{corollary}

Our next result gives a geometric way to construct conformal genuine
pairs by means of isometric ones, explaining the similitude between
\tref{int} and its isometric counterpart in \cite{df1}. To do this, we
need to introduce some further terminology.

Let $\Le^{N+2}$ denote the $(N+2)$--dimensional Lorentz space, and
let $$ \Ve^{N+1}= \{x\in\Le^{N+2}\colon\<x,x\>=0\} $$ be the light cone
in $\Le^{N+2}$. Fix a pseudo-orthonormal basis
$\{e_0,e_1,\ldots,e_{N+1}\}$ of $\Le^{N+2}$, that is,
$$
\<e_0,e_0\>=\<e_1,e_1\>=0,\ \<e_0,e_1\>=1
$$
and $\{e_2,\ldots, e_{N+1}\} $ is an orthonormal basis of the Riemannian
subspace $\{e_0,e_1\}^\perp$.
Then
$$
\Ee^{N}=\{x\in\Ve^{N+1}\colon\<x,e_0\>=1\}
$$
is a model of $N$--dimensional Euclidean space: the map
$\Psi\colon\;\R^{N}\to \Le^{N+2}$ defined by
\begin{equation}\label{e:psi}
\Psi(x) = -\frac{\|x\|^2}{2}e_0+ e_1+\sum_{i=1}^{N}x_ie_{i+1}
\end{equation}
is an isometric embedding with $\Psi(\R^{N})=\Ee^{N}$.

Given an immersion $g\colon\;M^n\to \Ve^{N+1}$ of a
differentiable manifold $M^n$, for any positive $\mu\in C^\infty(M^n)$ the
map $h\colon\;M^n\to \Ve^{N+1}$ given by $ h=\mu g$ is also an
immersion, and the induced metrics $\<\,,\,\>_g$ and $\<\,,\,\>_h$
are related by $\<\,,\,\>_h=\mu^2\<\,,\,\>_g$. Therefore, if
$\<\,,\,\>_g=\va^2\<\,,\,\>$ for some fixed metric $\<\,,\,\>$ on
$M^n$, then $h$ can be made isometric with respect to $\<\,,\,\>$
by choosing $\mu=\va^{-1}$. In particular, if $g=\Psi\circ f$ for
some conformal immersion $f\colon\,M^n\to \R^{N}$ of a
Riemannian manifold, then such an $h$ is denoted by ${\cal I}(f)$
and called the {\it isometric light cone representative\/} of $f$.

On the other hand, if $g\colon\;M^n\to \Ve^{N+1}$ is such that
$g(M^n)\subset \Ve^{N+1}\setminus\R_{e_0}$, where
$\R_{e_0}=\{t{e_0}:t>0\}$, define ${\cal C}(g)\colon
M^n\to{\mathbb {R}}^{N}$ by $\Psi\circ {\cal C}(g)
=\<g,e_0\>^{-1} g$. Since $\Psi$ is an isometric immersion, it
follows that $g$ and ${\cal C}(g)$ induce conformal metrics on
$M^n$ with conformal factor $\<g,e_0\>^{-1}$. Clearly, $g={\cal I}(f)$ if
$f={\cal C}(g)$ for an isometric immersion $g\colon\;M^n\to \Ve^{N+1}$ of
a Riemannian manifold.

This leads to the following procedure to construct a conformal pair of
immersions $f\colon\, M^n\to\R^{n+p}$ and
$\bar f\colon\,M^n\to\R^{n+q}$: start with a {\em Riemannian\/} manifold
$N^{n+1}$ that admits an isometric immersion
$F'\colon\,N^{n+1}\to\R^{n+p}$ and an isometric embedding
$\hat F\colon\,N^{n+1}\to\Le^{n+q+2}$ transversal to the light cone
$\Ve^{n+q+1}$. Then set
$M^n:=\hat F^{-1}(\hat F(N^{n+1})\cap \Ve^{n+q+2})$, $f=F'\circ i$ and
$\bar {f}={\cal C}(\hat F\circ i)$, where $i\colon\,M^n\to N^{n+1}$ is
the inclusion map.

The following result states that any genuine conformal pair $\{f,\bar f\}$
of Euclidean submanifolds in sufficiently low codimensions is locally
produced in this way from a genuine isometric pair $\{F,\hat F\}$ as
above.

\begin{theorem}\label{int2}\po
Assume that $f\colon M^n\!\to\R^{n+p}$, $p\geq 1$, and $\bar
f\colon M^n\to\R^{n+q}$ form a genuine conformal pair, with
$p+q\leq n-3$ and $\min\,\{p,q\}\le5$. Suppose further that $\bar
f$ is nowhere conformally congruent to an immersion that is
isometric to $f$. Then (locally on an open dense subset of $M^n$)
there exist a Riemannian manifold $N^{n+1}$ that admits an
isometric immersion $F'\colon\,N^{n+1}\to \R^{n+p}$ and an
isometric embedding $\hat F\colon\,N^{n+1}\to \Le^{n+q+2}$
transversal to the light cone $ \Ve^{n+q+1}$, and a conformal
diffeomorphism $i\colon\, M^n\to \hat F^{-1}(\hat F(N^{n+1})\cap
\Ve^{n+q+1})$ such that  $\{F',\hat F\}$ is a genuine isometric
pair,  $f=F'\circ i$ and $\bar {f}={\cal C}(\hat F\circ i)$.
\end{theorem}

Notice that the assumption that $\bar f$ is nowhere locally conformally
congruent to an immersion that is isometric to $f$ is always satisfied if
$f$ is genuinely {\it isometrically\/} rigid in $\R^{n+q}$, for instance
if $M^n$ does not carry any ruled open subset with rulings of dimension
at least $n-p-q$. In particular, this is always the case after composing
$f$ with a suitable inversion of $\R^{n+p}$.

\vskip .2cm

For $p=1$, Theorem \ref{int2} says that any hypersurface $f\colon
M^n\!\to\R^{n+1}$ that admits a genuine conformal (but not
isometric) deformation in $\R^{n+q}$ can be locally produced as
the intersection of an $(n+1)$-dimensional flat submanifold of
$\Le^{n+q+2}$ with the light cone:

\begin{corollary}\label{p1}\po
Let $f\colon M^n\!\to\R^{n+1}$ and $\bar f\colon M^n\to\R^{n+q}$ form a
conformal pair, with \mbox{$q\leq n-4$}. Assume that there exists no open
subset $M^n$ along which $\bar f$ is either a composition or it is
conformally congruent to an isometric deformation of $f$. Then, (locally
on an open dense subset of $M^n$) there exist an isometric embedding
$\bar F\colon\,U\subset \R^{n+1}\to\Le^{n+q+2}$ transversal to the light
cone $\Ve^{n+q+1}$ and a conformal diffeomorphism
$\tau\colon\,M^n\to \bar{M}^n:=\bar F^{-1}(\bar F(U)\cap \Ve^{n+q+1})
\subset U$ such that $f=i\circ \tau$ and
$\bar f={\cal C}(\bar F\circ \tau)$, where $i\colon\,\bar{M}^n\to U$ is
the inclusion map.
\end{corollary}

In the particular case $q=1$, the above reduces to Theorem $1$ in
\cite{dt1}, which can be regarded as a nonparametric description of
Cartan's conformally deformable hypersurfaces.

\vskip .2cm

Another important special case of \tref{int2} occurs when $q=0$. In this situation,
we consider a conformally flat submanifold
$f\colon\,M^n\to\R^{n+p}$, which clearly forms a genuine conformal pair
with any conformal diffeomorphism $\bar f\colon\,M^n\to U\subset \R^n$
onto an open subset. Then
we recover Theorem~1 from \cite{df4}, which gives a geometric
procedure to construct all conformally flat Euclidean submanifolds in
low codimension:

\begin{corollary}\label{cflat}\po
Let $f\colon M^n\!\to\R^{n+p}$, $n\geq 5$, $p\leq n-3$, be a conformal
immersion of a conformally flat manifold. Assume that the metric induced
by $f$ is nowhere flat. Then (locally on an open dense subset) there
exist a Riemannian manifold $N^{n+1}$ that admits an isometric immersion
$F\colon\,N^{n+1}\to \R^{n+p}$ and an isometric embedding
$\bar F\colon\,N^{n+1}\to \Le^{n+2}$, and a conformal diffeomorphism
$\tau\colon\,M^n\to \bar{M}^n:=\bar F^{-1}(\bar F(N^{n+1})\cap
\Ve^{n+1})$ such that $f=F\circ \tau$.
\end{corollary}

Our approach to study the geometric structure of a conformal pair of
immersions $f\colon M^n\!\to\R^{n+p}$ and $\bar f\colon M^n\!\to\R^{n+q}$
is, as usual, to fix on $M^n$ the Riemannian metric induced by one of the
immersions, say, $f$, and to reduce the problem to the study of the {\it
isometric} pair of immersions $\{f,\hat f\}$ that arises by considering
the isometric light cone representative $\hat f$ of $\bar f$. However,
two distinct cases need to be handled separately, and a delicate
degenerate case requires to deal with the isometric light cone
representatives of both $f$ and $\bar{f}$. Therefore, as a first and main
step, in the next section we extend the theory developed in \cite{df1} to
isometric pairs of immersions into flat {\it semi}-Riemannian spaces.

\vskip .5cm

Before we conclude this introduction, one final remark is in order.
Although we deal only with pairs of immersions in this paper, when
studying some rigidity phenomena of submanifolds one is naturally led to
consider {\it sets} of immersions, not only pairs. For instance, it was
shown in Theorem 5 in \cite{df2} that the associated family
$\{f_\theta:M^{2n}\to\R^{2n+2}: \theta\in[0,\pi)\}$ of a minimal
nonholomorphic Kahler submanifold $f_0$ of rank two in codimension two
does not extend isometrically, although for any $\theta_1\neq\theta_2\in
[0,\pi)$ the pair $\{f_{\theta_1},f_{\theta_2}\}$ does extend. Here, we
say that a set $\{f_i:M^n \to N_i^{n+p_i}: i\in I\}$ of isometric (resp.,
conformal) immersions indexed by an arbitrary set $I$ {\it extends
isometrically} (resp., {\it conformally}), when there exist an isometric
(resp., conformal) embedding $j: M^n \to N^m$, with $m > n$, and a set
$\{F_i:N^m \to N_i^{n+p_i}: i\in I\}$ of isometric (resp., conformal)
immersions such that $f_i=F_i\circ j$, for all $i\in I$.

\section{Isometric pairs into flat semi-Riemannian spaces}  

In this section, we study the structure of the tangent and normal bundles
of a pair of isometric immersions into flat semi-Riemannian spaces. Our
goal is to give conditions that allow the construction of isometric ruled
extensions.

\subsection{Semi-Riemannian ruled isometric extensions}  
In this subsection we give general conditions for the existence of
isometric ruled extensions of a pair of isometric immersions into flat semi-Riemannian spaces.
The proofs  are identical to the ones for the Riemannian case
(\cite{df1}) and will be omitted.

\vskip .5cm

Throughout the paper, given a bilinear form $\beta\colon\, V^n\times
U^m\to W$ between finite dimensional real vector spaces, we denote by
$\Sal(\beta)\subset W$ the subspace spanned by the image of~$\beta$,
that is,
$$
\Sal(\beta)=\spa\{\beta(X,Y):\, X\in V^n\;
\mbox{and}\; Y\in U^m\},
$$
and by $\N(\beta)\subset V^n$ the (left) nullity space of $\beta$ defined
as
$$
\N(\beta) =\{X\in V^n :\,\beta(X,Y)=0\; \mbox{for all}\,\; Y\in U^m\}.
$$
If $W$ is endowed with a nondegenerate
inner product $\<\ ,\ \>$ and $T\subseteq W$ is a nondegenerate
subspace with respect to $\<\ ,\ \>$, we denote $\beta_T=\pi_T\circ\beta$,
where $\pi_T$ is the orthogonal projection onto $T$. Then
$$
\N(\beta_T)=\{X\in V^n :\,\<\beta(X,Y),\xi\>=0\;
\mbox{for all}\,\; Y\in U^m, \xi\in T\}.
$$
We also denote by $\N(\beta_T)$ the subspace defined as above even if $T$
is degenerate.

Let $\R_a^m$ stand for $\R^m$ with the standard flat
semi-Riemannian metric of index $a$. In particular,
$\R_1^m=\Le^m$. Let $f\colon\,M^n\to\R_a^{n+p}$ and
$\hat{f}\colon\,M^n\to\R_b^{n+q}$ form an isometric pair. Assume
that there exists a vector bundle isometry
$$
\Tes\colon\,L\subset T^{\perp}_f M\to\hat{L}\subset T^{\perp}_{\hat f}M
$$
  between nondegenerate subbundles such that
$$
D=\N(\a_{L^\perp})\cap \N(\hat\a_{\hat{L}^\perp})\subset TM
$$
defines a smooth subbundle of $TM$ and such that the pair $(\Tes,D)$ satisfies the
following two conditions:
\begin{equation}\label{ces}
\left\{
\begin{array}{l}\!(C_1)\;  \mbox{The isometry}
\;\Tes\; \mbox{is parallel and preserves second fundamental forms;}
\vspace*{2ex}\\
\!(C_2)\;  \mbox{The subbundles}\; L \;\mbox{and} \;\hat L\;
\mbox{are parallel along}\; D \;\mbox{in the  normal connections.}
\end{array} \right.
\end{equation}

\vspace {1ex}

Let $\phi\colon (TM\oplus L)\times TM\to L^\perp\times \hat L^\perp$
be the bilinear form given by
\begin{equation}\label{bili}
\phi(Y+\xi,X) =
\left((\nab_X(Y+\xi))_{L^\perp},(\nab_X(Y+\Tes\xi))_{\hat L^\perp}\right),
\end{equation}
where $\nab$ stands for the connections of both $\R_a^{n+p}$ and
$\R_b^{n+q}$, and assume further that the vector subspaces
$$
\Delta:=\N(\phi)\subset TM\oplus L
$$
have constant dimension on $M^n$. By condition $(C_1)$, the vector bundle
isometry defined as
$T_0=I\oplus \Tes\colon\,f_*TM\oplus L \to\hat{f}_*TM\oplus\hat{L}$
is parallel with respect to the connections induced by the Euclidean
ambient spaces. It follows that
$T_0\,|_\Delta\colon\,\Delta\to\hat\Delta$ is a parallel vector bundle
isometry, and hence, we may identify $\hat\Delta$~with~$\Delta$.

\begin{lemma}\label{inter}{\rm (\cite{df1}).} The distribution
$D\subset\Delta$ is integrable and $\Delta\cap TM=D$ holds.
\end{lemma}

\vspace{1ex}

Consider the vector bundle $\pi\colon\,\Lambda=\hat\Lambda\to M^n$
determined by the orthogonal splitting $\Delta = D\oplus\Lambda$,
and define $F'\colon\, N\to\R_a^{n+p}$ as the restriction of the
map
\begin{equation}\label{e:par}
\lambda\in\Lambda\mapsto f(\pi(\lambda))+ \lambda
\end{equation}
to a tubular neighborhood $N$ of the $0$-section
$j\colon\,M^n\hookrightarrow N\subset\Lambda$ of $\Lambda$ along
which $F$ is an immersion. Similarly, define $\hat F\colon\,
N\to\R_b^{n+q}$. Henceforth, $L^\perp$ and $\Delta$ will be
understood as vector bundles over $N\subset\Lambda$ by means of
$L^\perp(\lambda)=L^\perp (\pi(\lambda))$ and
$\Delta(\lambda)=\Delta(\pi(\lambda))$.

\begin{proposition}\label{FF}{\rm (\cite{df1}).} The immersions
$F'$ and $\hat F$ are isometric \mbox{$\Delta$-ruled} extensions
of $f$ and~$\hat{f}$. Moreover, there are  smooth orthogonal
splittings
\begin{equation}\label{split}
T_{F'}^\perp N=\Les\oplus L^\perp\;\;\;\mbox{and}\;\;\;
T_{\hat F}^\perp N=\hat \Les\oplus \hat{L}^\perp
\end{equation}
and a vector bundle isometry $T\colon\,\Les\to\hat{\Les}$
such that
\begin{equation}\label{inc}
\Delta=
\N(\a^{F'}_{\Les^\perp})\cap \N(\hat\a^{\hat F}_{\hat{\Les}^\perp}),
\end{equation}
and the pair $(T,\Delta)$ satisfies conditions $(C_1)$ and
$(C_2)$ in {\em(\ref{ces})}.
\end{proposition}

Observe that if the ruled extensions $F'$ and $\hat F$ are
{\it trivial\/} (i.e., $\dim N = n$) then $f$ and $\hat f$
are themselves $D$-ruled.

\subsection{Construction of the pair $(\Tes,D^d)$ and the estimate on $d$}  

In this subsection we show how to construct a pair $(\Tes,D^d)$
satisfying conditions $(C_1)$ and $(C_2)$ in {(\ref{ces})} for a pair of
isometric immersions $f\colon\,M^n\to\R_a^{n+p}$ and
$\hat{f}\colon\,M^n\to\R_b^{n+q}$, and we obtain an estimate on $d$. We
follow closely the strategy in \cite{df1} for the case $a=0=b$. Here,
however, two distinct cases arise, depending on whether a certain
nondegeneracy condition is satisfied or not. The degenerate case requires
several modifications in the arguments of \cite{df1}, which will be
carried out in Subsection $2.2.2$ only in the case that is needed for our
study of conformal pairs of immersions into Euclidean space.

\vskip .3cm

Given an isometric pair of immersions $f\colon\,M^n\to\R_a^{n+p}$ and
$\hat{f}\colon\,M^n\to\R_b^{n+q}$, denote by $\a$ and $\hat \a$ their
respective second fundamental forms and endow the vector bundle
$T_f^\perp M\oplus T_{\hat{f}}^\perp M$ with the indefinite metric of
type $(p,q)$ given by
$$
\<\!\<\;\,,\;\,\>\!\>_{T_f^\perp M\oplus T_{\hat{f}}^\perp M}=
\<\,\;,\;\,\>_{T_f^\perp M}
- \<\,\;,\;\,\>_{T_{\hat{f}}^\perp M}.
$$
Set $\a\oplus\hat\a\colon\,TM\times TM\to\Sal(\a)
\oplus\Sal(\hat\a)\subset T_f^\perp M \oplus T_{\hat{f}}^\perp M$.
\begin{definition}\po\label{d:nd}
{\rm
We say that the pair $\{f,\hat f\}$ as above is {\it nondegenerate}
if the projections of
$\Omega=\Omega(f,\hat f):=\Sal(\a\oplus\hat\a)\cap
\Sal (\a\oplus\hat\a)^\perp\subset\Sal(\a)\oplus\Sal(\hat\a)$
onto $T_f^\perp M$ and $T_{\hat f}^\perp M$
are injective. When this condition is nowhere satisfied,
we say that the pair is {\it degenerate}.
}
\end{definition}
Notice that $\{f,\hat f\}$ is nondegenerate if both
$\Sal (\a)$ and $\Sal (\hat\a)$ are nondegenerate. In particular,
this is the case if  $a$, $b$ and the index of $M^n$ are all equal.

\subsubsection{The nondegenerate case}

Assuming $\{f,\hat f\}$ to be nondegenerate, we have orthogonal
splittings
$$
\Sal(\a)=\Gamma\oplus\Gamma^\perp\;\;\;\mbox{and}
\;\;\;\;\Sal(\hat\a)=\hat{\Gamma}\oplus
\hat{\Gamma}^\perp,
$$
where $\Gamma=\Sal(\a)\cap\Omega^\perp$ and
$\hat\Gamma=\Sal(\hat\a)\cap\Omega^\perp$, and an
isometry
$\Ies\colon\,\Gamma^\perp\to\hat{\Gamma}^\perp$ such that
$$
\Omega=\{(\eta,\Ies\eta):\eta\in\Gamma^\perp\}\subset \Gamma^\perp\oplus\hat
\Gamma^\perp
$$
and $\hat\a_{\hat{\Gamma}^{\perp}}
=\Ies\circ\a_{\Gamma^{\perp}}$.
 From now on we identify $\Gamma^\perp$ with
$\hat{\Gamma}^\perp$ by means of $\Ies$, and hence
\begin{equation}\label{*}
\hat\a_{\hat{\Gamma}^{\perp}}=\a_{\Gamma^{\perp}}.
\end{equation}

Define $\beta\colon\,TM\times TM\to\Gamma\oplus\hat\Gamma $
as $\beta=\a_{\Gamma}\oplus\hat\a_{\hat\Gamma}$,
and a vector subbundle $\Theta\subset TM$ by
$$
\Theta=\N(\beta).
$$
The vector subbundle $S\subset\Gamma^\perp(=\hat\Gamma^\perp)$ defined by
\begin{equation}\label{e:s}
S=\Sal(\a|_{\Theta\times TM})
\end{equation}
satisfies $\Theta=\N(\a_{S^\perp})\cap\N(\hat\a_{\hat{S}^\perp})$.
Now, given $X\in TM$, denote by $\Kes(X)\in\Lambda^2(S)$  the
skew-symmetric tensor given by
$$
\Kes(X)\eta=(\nabla^\perp_X\eta)_{S}-
(\hat{\nabla}^{\perp}_X\eta)_{\hat S},
$$
and  define a vector subbundle $S_0\subset S$ by
$$
S_0=\bigcap_{X\in TM}\ker \Kes(X).
$$
Then, define vector subbundles $L^\ell\subset S_0$ and $D^d\subset\Theta$
as
\begin{equation}\label{e:ell}
L^\ell=\{\delta\in S_0:\, \nap_Y\delta\in S\;\;\mbox{and}\;\;
\hat{\nabla}^\perp_Y\delta\in \hat S\;\;\mbox{for all}\;\; Y\in\Theta\}
\end{equation}
and
$$
D^d=\N(\a_{L^\perp}) \cap\N(\hat\a_{\hat{L}^\perp}),
$$
and let $\Tes\colon\,L^\ell\to L^\ell$ be the induced vector bundle
isometry given by
$$
\Tes=\Ies|_L\colon\,L^\ell\subset T^{\perp}_f M\to
L^\ell\subset T^{\perp}_{\hat f}M.
$$

With these definitions, we have the following:

\begin{theorem}\label{t3}{\rm (\cite{df1}).}
Let $f\colon M^n\!\to\R_a^{n+p}$ and $\hat{f}\colon M^n\to\R_b^{n+q}$
form a nondegenerate isometric pair of immersions of a semi-Riemannian
manifold with index $\min\{a,b\}$. Then, along each connected component
of an open dense subset of $M^n$, the pair $(\Tes, D^d)$ satisfies
$(C_1)$ and $(C_2)$ in {\em(\ref{ces})}.
In particular, $f$ and $\hat f$ have {\em (}possibly trivial\/{\em)}
maximal isometric $\Delta^{d+r}$-ruled extensions
$F'\colon\,N^{n+r}\to\R_a^{n+p}$ and $\hat F\colon\,N^{n+r}\to\R_b^{n+q}$,
$0\leq r\leq\ell$, which form a nondegenerate pair and satisfy the
conclusions of \pref{FF}.

Moreover, if $p+q\leq n-1$ and
$\min\,\{p+b-a,q+a-b\}\le 6$ then
\begin{equation}\label{e:gen}
d+r\geq n-p-q+3\ell,
\end{equation}
unless $\min\,\{p+b-a,q+a-b\}=6$ and $\ell=0$ in which case
$d+r\geq n-p-q+3\ell-1$.
\end{theorem}

\begin{remark}\po\label{r:s} {\rm
$i)$ The hypothesis on the codimensions in \tref{t3}
is required in a fundamental algebraic result needed in its proof,
whose most general version is Theorem~3 in \cite{df3}.
Unfortunately, this algebraic result is false without that assumption
(\cite{df5}).

$ii)$ We can relax the hypothesis on the index of $M$
by asking $S$ in \eqref{e:s} to be Riemannian.
}
\end{remark}

The proof of the above result follows exactly as those of Theorems 11 and
14 in \cite{df1}, where of course no hypothesis on the nondegeneracy was
needed since both normal spaces are Riemannian when both ambient spaces
are Euclidean. The only extra property to verify is that the pair
$\{F',\hat F\}$ in \tref{t3} is nondegenerate, but this is immediate from
\pref{FF}, since
$\a^f_{L}=\a^{\hat f}_{\hat L}$,
$\a^{F'}_{\Les} = \hat\a^{\hat F}_{\hat\Les}$,
$\a^{F'}_{\Les^\perp}|_{TM\times TM} = \a^f_{L^\perp}$,
$\a^{\hat F}_{\hat\Les^\perp}|_{TM\times TM}=\a^{\hat f}_{\hat L^\perp}$
and
$\Delta=\N(\a^{F'}_{\Les^\perp})\cap\N(\hat\a^{\hat F}_{\hat\Les^\perp})$.

\vspace{2ex}

\subsubsection{The degenerate case}

In this subsection we address the degenerate case in the setting that
will be needed for the next section, namely, for a pair of isometric
immersions $f\colon\,M^n\to\R^{n+p}$ and
$\hat{f}\colon\,M^n\to\Ve^{n+q+1}\subset\Le^{n+q+2}$ of a Riemannian
manifold, under the assumptions that
\begin{equation}\label{e:dimsc}
p+q\leq n-1 \ \ \ {\rm and} \ \ \ \min\,\{p,q\}\le 5.
\end{equation}
The main difficulty here is that the subbundle $\Gamma^\perp$,
constructed in the preceding subsection for the nondegenerate case, is no
longer well defined. To deal with this issue, we need to consider also
the isometric representative of $f$ into the light cone.

\vskip .3cm

Thus, here we assume that
there is a null vector field
$0\neq \xi_0\in \Sal(\hat \a)\cap \Sal(\hat \a)^\perp$
such that
$(0,\xi_0) \in \Sal (\a\oplus\hat\a)\cap\Sal (\a\oplus\hat\a)^\perp$.
Set $f'={\cal I}(f)\colon\,M^n\to\Ve^{n+p+1}\subset\Le^{n+p+2}$ and
denote by $\a'$ its second fundamental form. Then, the position vector
fields of both $f'$ and $\hat f$ are normal, and
\begin{equation}\label{e:ssf}
A^{f'}_{f'} = A^{\hat f}_{\hat f} = -I,\ \ \
A^{f'}_{e_0} = A^{\hat f}_{\xi_0} = 0.
\end{equation}
Observe that $(f',\hat f)\in\Sal(\a'\oplus\hat\a)^\perp$, but $(f',\hat
f)\not\in\Sal(\a'\oplus\hat\a)$. Since the normal spaces of $\hat f$ have index $1$,
we can assume further that $\<\hat f,\xi_0\>=1$.

\vspace {1ex}

We will make a similar construction as in the nondegenerate case, but now
for the pair $\{f',\hat f\}$. The idea is to force the inclusion of the
vector $(f',\hat f)$ in $\Omega$, despite the fact that
$(f',\hat f)\not\in \Sal(\a'\oplus\hat\a)$. We then define
$\Omega\subset\spa\{(f',\hat f)\}\oplus\Sal(\a')\oplus\Sal(\hat\a)$
as the vector bundle with null fibers
$$
\Omega:=\spa\{(f',\hat f)\}\oplus\left(\Sal
(\a'\oplus\hat\a)\cap\Sal(\a'\oplus\hat\a)^\perp\right).
$$
Notice that, by \eqref{e:ssf} and the definition of $\xi_0$, we
get that $(e_0,\xi_0)\in\Omega$. As before, there are orthogonal
splittings
$$
\spa\{f'\}\oplus\Sal(\a')=\Gamma\oplus\Gamma^\perp\;\;\;\mbox{and}
\;\;\;\;\spa\{\hat f\}\oplus\Sal(\hat\a)=\hat{\Gamma}\oplus
\hat{\Gamma}^\perp,
$$
where $\Gamma=\Sal(\a')\cap\Omega^\perp\subseteq\spa\{f',e_0\}^\perp$ and
$\hat\Gamma=\Sal(\hat\a)\cap\Omega^\perp\subseteq\spa\{\hat f,\xi_0\}^\perp$
are Riemannian, and an isometry
$\Ies\colon\,\Gamma^\perp\to\hat{\Gamma}^\perp$ such that
$$
\Omega=\{(\eta,\Ies\eta):\eta\in\Gamma^\perp\}\subset
\Gamma^\perp\oplus\hat \Gamma^\perp,
$$
with $\hat\a_{\hat{\Gamma}^{\perp}}=\Ies\circ\a'_{\Gamma^{\perp}}$
giving the same identification as before. Notice that
$\Ies(f')=\hat f$ and $\Ies(e_0) = \xi_0$.

The preceding {\it ad-hoc} inclusion in $\Gamma^\perp$ and
$\hat\Gamma^\perp$ of the position vectors $f'$ and $\hat f$,
respectively, despite the fact that they are not contained in the
subspaces spanned by the images of $\a'$ and $\hat \a$, requires a few
arguments to show that some properties of the bundles used in the proofs
in \cite{df1}, which were automatic in the Riemannian case, still hold in
our situation. We will prove these in the form of numbered claims.

Define $\beta$ and $\Theta$ as in the nondegenerate case, but
$S\subseteq\Gamma^\perp(=\hat\Gamma^\perp)$ as
$$
S=\spa\{f'\}\oplus\Sal(\a'|_{\Theta\times TM}).
$$

\vspace{2ex}
\noindent {\it Claim 1. The subbundle $S$ is Lorentzian.}
\vspace{1ex}

\proof
As in the beginning of the proof of Lemma 13 in \cite{df1}, we easily
check that $\Sal(\beta)$ is nondegenerate. By \eqref{e:dimsc}
and Theorem 3 in \cite{df3} (see also Corollary 17 in \cite{df1}), we have
that $\dim \Theta \geq n-\dim \Sal(\beta) \geq n-p-q > 0$. The claim
follows from the definition of $S$ and the fact that $f'\in S$ is null,
since by \eqref{e:ssf},
\begin{equation}\label{e:th0}
\<\a'(Z,Z),f'\> = -\|Z\|^2 \neq 0, \fall 0\neq Z\in\Theta.\qed
\end{equation}

We also have that
$\Theta=\N(\a'_{S^\perp})\cap\N(\hat\a_{\hat{S}^\perp})$.
Then, we define  $S_0,\Kes,L^\ell,\Tes$ and $D^d$ just as before.
Observe that, since $f'$ and $\hat f$ are normal parallel, we obtain
\begin{equation}\label{e:fl}
f'\in L^\ell \subseteq S_0.
\end{equation}

\vspace{2ex}
\noindent {\it Claim 2. The subbundle $S_0\subseteq S$ is Lorentzian.}
\vspace{1ex}

\proof
In our setting, the proof of Lemma 12 in \cite{df1} implies  that the tensor $\Kes$ as a map
\mbox{$\Kes\colon\, TM\to\Lambda^2(S)$} satisfies that
$\Ima\Kes(Z)\subset S\cap\Sal(\a'|_{\Theta\times TM})^\perp$,
for all $Z\in\Theta$. But since \mbox{$f'\in\kerl\Kes(Z)=(\Ima\Kes(Z))^\perp$},
we obtain that $\Ima\Kes(Z)\subset S\cap S^\perp = 0$, that is,
$$
\Kes(Z)=0, \fall Z\in\Theta,
$$
and the statement of Lemma 12 in \cite{df1} holds here also. But then $\a'(Z,Z)\in S_0$,
and the result follows from \eqref{e:th0} and the fact that $f'\in S_0$
is null.
\qed
\vspace{1.5ex}

By the above two claims and the antisymmetry of $\Kes$, exactly as in
\cite{df1} we have
the orthogonal splitting
$$
S=S_0\oplus^\perp S_1,
$$
with $S_1=\spa\{\Kes(X)S_1: X\in TM\}$, which is a Riemannian subbundle. Moreover,
by Claim 2 and \eqref{e:dimsc} we have that
$\dim S_1\leq 5$.

\vspace{2ex}
\noindent {\it Claim 3. The bilinear map
$\gamma=\a'_{S_1}|_{\Theta\times TM}$ satisfies $S_1=\Sal(\gamma)$.}

\vspace{1ex}
\proof
Consider $\xi\in S_1\cap\Sal(\gamma)^\perp$. Hence,
$\xi\in \Sal(\a'|_{\Theta\times TM})^\perp$ since $S_0\perp S_1$.
But then $\xi\perp f'\in S_0$, and $\xi\in S\cap S^\perp=0$.
\qed
\vspace{1.5ex}

\vspace{2ex}
\noindent {\it Claim 4. It holds that $d>0$,
and the subbundle $L^\ell$ is Lorentzian.}
\vspace{1ex}

\proof Using the previous claims, one can argue exactly as in the
proof of Theorem 14 in \cite{df1} to conclude that
$$
d\geq n-p-q + 2\ell>0,
$$
bearing in mind  Remark 20 part 3 in
\cite{df1}. Now, we get from the definition of $D$ and \eqref{e:th0}
applied to $0\neq Z\in D$ that $f'\not\in L^\perp$. The result now
follows from the fact that $f'\in L$.
\qed
\vspace{1.5ex}

These claims are all we need to make a straightforward check that the
proofs of Theorems 11 and 14 in \cite{df1} still work in the degenerate
case with the preceding definitions.
We obtain that, along each connected component of an open dense subset of
$M^n$, the pair $(\Tes, D^d)$ for $\{f',\hat f\}$ satisfies $(C_1)$ and
$(C_2)$ in \eqref{ces}, and that the immersions $f'$ and $\hat f$ have
mutually $\Delta^{d+r}$-ruled isometric extensions $F'\colon
N^{n+r}\to\Le^{n+p+2}$ and $\hat F\colon\, N^{n+r}\to\Le^{n+q+2}$. Since
$f'\in L^\ell$, we have that $f'\in\Delta^{d+r}=\N(\phi)$ as in
\eqref{bili}. Thus, by \eqref{e:par}, both $F'(N)$ and $\hat F(N)$ are
cones, where a subset ${\cal C}\subseteq \R^m$ being a {\it cone} means
that $x\in {\cal C}$ implies $tx\in {\cal C}$ for $t$ close to~1.

Observe also that the extensions $F'$ and  $\hat F$ are Lorentzian, since
$f'\in\Delta^{d+r}\subset TN$ and, for $0\neq Z\in D^d$, we get that
$\tilde\nabla_ZZ\in\Delta^{d+r}$ and $\<\tilde\nabla_ZZ,f'\> =
\<\a'(Z,Z),f'\>=-\|Z\|^2<0$. We also conclude that $\Delta^{d+r}$ is
strictly larger than $D^d\oplus\spa\{f'\}$ (and is, in fact, a Lorentzian subbundle),
hence \lref{inter} implies that $r\geq 2$.

Moreover, by the observation after \dref{d:nd}, the pair $\{F',\hat F\}$
is nondegenerate. Hence, under the codimension assumption
\eqref{e:dimsc}, we may apply \tref{t3} to $\{F',\hat F\}$ to conclude
that
$$d+r \geq n+r - (p+2-r)- (q+2-r) + 3(\ell-r) = n-p-q+3\ell-4.$$

 Summarizing, we have the following result:

\begin{proposition}\label{deg1}\po
Let $f\colon M^n\!\to\R^{n+p}$ and
$\hat{f}\colon M^n\to\Ve^{n+q+1}\subset\Le^{n+q+2}$ form a degenerate
isometric pair of immersions, and set
$f'={\cal I}(f)\colon\, M^n\to\Ve^{n+p+1}\subset\Le^{n+p+2}$. Assume that
$p+q\leq n-1$ and $\min\,\{p,q\}\le 5$. Then, along each connected
component of an open dense subset of $M^n$, the immersions $f'$ and $\hat f$ have
mutually $\Delta^{s}$-ruled isometric Lorentzian conical extensions
$F'\colon N^{n+r}\to\Le^{n+p+2}$ and
$\hat F\colon\,N^{n+r}\to\Le^{n+q+2}$ such that
$\<F',F'\>=\<\hat F,\hat F\>$, with
$$
s\geq n-p-q+3\ell-4, \ \ 2\leq r\leq \ell.
$$
Moreover, there exists a vector bundle isometry
$T:\Les^{l-r}\to\hat\Les^{l-r}$
satisfying the conclusions of Proposition \ref{FF}.
\end{proposition}

\begin{remark}\po\label{r:vale}
{\rm \tref{t3} still holds if the pair is degenerate without the
codimension assumption \eqref{e:dimsc}, but the proof is not completely
analogous to the one in \cite{df1}. We omit it here since it is not needed
for our purposes.}
\end{remark}

\section{The main result}

In this section we will prove a slightly more general version of
\tref{int}. We start with the following preliminary fact.

\begin{lemma}\po\label{l:sff}
Let $f\colon\,M^n \to \R^{n+p}$ be a conformal immersion of a Riemannian
manifold, and let
$f'={\cal I}(f)\colon\, M^n \to \Ve^{n+p+1}\subset \Le^{n+p+2}$ be its
isometric light cone representative. If $f'$ is $\Delta$-conformally
ruled, then the same holds for $f$. Moreover, the following holds:
\begin{itemize}
\item[$(i)$] There exists $\lambda\in C^{\infty}(M)$ such that the
conformal factor $\va\in C^{\infty}(M)$ relating the metrics
$\<\,,\,\>_{f}$ and $\<\,,\,\>_{f'}$ satisfies
$\hess\va|_{\Delta\times\Delta}=\lambda\<\,,\,\>_{f'}|_{\Delta\times\Delta}$;
\item[$(ii)$] The normal components of the mean curvature vector fields
$\eta$ and $\eta'$ of the leaves of $\Delta$ for $f$ and $f'$ are related
by
\begin{equation}\label{e:etas}
\eta'=\va^{-1}(d\Psi(\eta)-\va\xi+\lambda {f'}),
\end{equation}
where $\xi=\va^{-1}e_0-d(\Psi\circ f)\ {\rm grad}\, \va$;
\item[$(iii)$] The symmetric bilinear forms
$\beta^f=\alpha^{{f}}-\<\,,\,\>_f{\eta}$ and
$\beta^{f'}=\alpha^{{f'}}-\<\,,\,\>_{f'}{\eta'}$ are related by
\begin{equation}\label{e:sffs3}
\beta^{f'}=
\va d\Psi(\beta^f)+\va^{-1}(\hess \va(\,\,,\,)-\lambda\<\,,\,\>_{f'}){f'},
\end{equation}
\end{itemize}
where the Hessian and the gradient are computed with respect to
$\<\,,\,\>_{f'}$.
\end{lemma}
\proof
By Lemma $4$ in \cite{to}, the second fundamental forms of ${f}$ and
${f'}$ are related by
\begin{equation}\label{e:sffs}
\alpha^{{f'}}(\,,\,)=
d\Psi(\va\alpha^{{f}}(\,,\,))-\<\,,\,\>_{f'}\xi+\va^{-1}\hess \va(\,,\,){f'}.
\end{equation}
Since ${f'}$ is $\Delta$-conformally ruled, there exists a  normal vector
field $\eta'\in T^\perp_{{f'}}M$ such that
$\alpha^{{f'}}|_{\Delta\times \Delta}=
\<\,,\,\>_{{f'}}|_{\Delta\times \Delta}{\eta'}$.
It follows that there exist a normal vector field
$\eta\in T^\perp_{{f}}M$ and $\lambda\in C^{\infty}(M)$ such that
$\alpha^{{f}}|_{\Delta\times \Delta}=
\<\,,\,\>_f|_{\Delta\times \Delta}{\eta}$ and
$\hess \va|_{\Delta\times \Delta}=
\lambda\<\,,\,\>_{f'}|_{\Delta\times \Delta}$.
Replacing into \eqref{e:sffs} yields \eqref{e:etas} and \eqref{e:sffs3}.
\qed

\begin{theorem}\label{int3}\po
Let $f\colon M^n\!\to\R^{n+p}$ and $\bar{f}\colon M^n\to\R^{n+q}$
form a conformal pair, with $p+q\leq n-3$ and $\min\,\{p,q\}\le
5$. Then (locally on an open dense subset of $M^n$) the pair
$\{f,\bar{f}\}$ extends conformally (possibly trivially) to a
mutually $\Delta^s$-conformally ruled pair of immersions $F\colon
N^{n+r}\!\to\R^{n+p}$ and $\bar F\colon N^{n+r}\!\to\R^{n+q}$,
with
$$
\Delta^s=
\N(\beta^{F}_{\Les^\perp})\cap \N(\beta^{\bar F}_{\bar\Les^\perp})
\ \ \ \ \ and \ \ \ \ \ s\geq n-p-q+3(\ell^c+r),
$$
where $\Les:=L^c_\Delta(F)$, $\bar \Les:=L^c_\Delta(\bar F)$
and $\ell^c:=\rank\Les=\rank\bar\Les$.
Moreover, there exists a parallel
vector bundle isometry $\Tes\colon\, \Les \to \bar \Les$ such that
$\Tes\circ \beta^{F}_{\Les}=\va\beta^{\bar{F}}_{\bar \Les}$,
where $\va$ is the conformal factor relating the metrics induced
by $\bar F$ and $F$.
\end{theorem}
\proof Set
$\hat{f}={\cal I}(\bar{f})\colon M^n\to\Ve^{n+q+1}\subset\Le^{n+q+2}$,
so that $\{f,\hat{f}\}$ becomes an isometric pair. We consider separately
the two possible cases:

\vspace{1ex}
\noindent {\it $i)$ The pair $\{f,\hat f\}$ is nondegenerate.}
\vskip .2cm

In this case,  Theorem \ref{t3} applies and yields maximal isometric
$\Delta_0^{s_0}$-ruled extensions $F'\colon\,N_0^{n+r_0}\to\R^{n+p}$ and
$\hat F\colon\,N_0^{n+r_0}\to\Le^{n+q+2}$,
which form a nondegenerate pair and satisfy the
conclusions of \pref{FF}.
Since $\hat{f}$ takes values in the light cone and $M^n$ is Riemannian,
$\hat{f}$ can not be ruled, and thus $r_0\geq 1$. In particular, we have
$\ell=\ell_0+r_0\geq 1$, and hence \eqref{e:gen} gives
$$
s_0\geq n-p-q-2+3(\ell_0+r_0),
$$
where $\ell_0$ is the rank of the subbundle $L$ given by (\ref{e:ell})
for the pair $\{F',\hat{F}\}$, which we denote by ${\Les}_0$.
Moreover, $\hat{F}$ is transversal to the light cone, for $N_0^{n+r_0}$
is Riemannian and $\hat{F}$ is ruled. By restricting to an open subset if
necessary, we may assume that $\hat{F}$ is an embedding, and hence
$N:=\hat{F}^{-1}(\hat{F}(N_0)\cap \Ve^{n+q+1})\supseteq M^n$ is an
$(n+r_0-1)$-dimensional manifold.

Set $F=F'\circ i$ and $\bar{F}= {\cal C}(\hat{F}\circ i)\colon\,{N}\to
\R^{n+q}$, where $i\colon\, N\to N_0$ is the inclusion map. Then
$\{F,\bar{F}\}$ is a conformal pair, $F\circ j=f$ and $\bar{F}\circ
j=\bar{f}$, where $j$ is the inclusion of $M$ into $N$, and hence
$\{F,\bar{F}\}$ is a conformal extension of $\{f,\bar f\}$. Moreover, $F$
and $\bar{F}$ are mutually $\Delta^s$-conformally ruled, where $s=s_0-1$
and $\Delta^s$ is the distribution on $M$ defined by
$d\hat{F}(\Delta)=d\hat{F}(\Delta_0)\cap T\Ve^{n+q+1}$.
Therefore,
$$
s\geq n-p-q+3(\ell_0+r),
$$
with $r=r_0-1$, and hence the estimate on $s$ will follow once we prove
that $\ell_0\geq \ell^c$.

First observe that
$\hat{\Les}:=L^c_{\Delta}(\hat{F}\circ i)\subset
\hat{\Les}_0\oplus \spa\{d\hat F(\eta)\}$,
where $\eta(x)\in T_xN_0\cap T_x^\perp N $ is the normal component of the
mean curvature vector at $x\in N$ of the leaf of $\Delta^s$ through~$x$.
On the other hand, we obtain from \eqref{e:sffs3} that
$d\Psi(\bar{\Les})\subset \hat{\Les}\oplus \spa\{{\hat F}\circ i\}$,
hence
\begin{equation}\label{e:lor} d\Psi(\bar{\Les})\subset
\hat{\Les}_0\oplus \spa\{d\hat F(\eta)\}\oplus \spa\{{\hat F}\circ i\}.
\end{equation}

Now, for a unit vector $Z\in\Delta$, we have
$\alpha_{\hat{F}\circ i}(Z,Z)=d\hat F(\eta)$.
So $\<d\hat F(\eta),\hat{F}\circ i\>=-1$, since
$\<\alpha_{\hat{F}\circ i}(\,\,,\,),\hat{F}\circ i\>=-\<\,\,,\,\>$.
It follows that the subspace $V$ on the right-hand-side of \eqref{e:lor}
is Lorentzian. Therefore, the subspace $d\Psi(\bar{\Les})$ is a
Riemannian subspace of $V$ that is orthogonal to the null vector ${\hat
F}\circ i\in V$, hence $d\Psi(\bar{\Les})$ has codimension at least
two in $V$. We conclude that $\rank \bar\Les \leq \ell_0$, as we wished.

We show now that
\begin{equation}\label{e:iso}
\Tes(\beta^{F}(Z,X))=\va\beta^{\bar F}(Z,X),\,\,\,\,\mbox{for}\,\,Z\in
\Delta,\,X\in TN,
\end{equation}
defines a parallel vector bundle isometry $\Tes\colon\, \Les\to\bar\Les$
such that $\Tes\circ\beta^{F}_{\Les}=\va\beta^{\bar{F}}_{\bar \Les}$.

Extend $T_0$ to a vector bundle map $T_1$ between
${\Les}_1:={\Les}_0\oplus \spa\{\eta_F\}\subset T_F^\perp N$ and
$\hat{\Les}_1:=\hat{\Les}_0\oplus \spa\{\eta_{\hat{F}\circ i}\}\subset
T_{\hat{F}\circ i}^\perp N$
by setting $T_1|_{{\Les}_0}=T_0$ and $T_1(\eta_F)=\eta_{\hat{F}\circ i}$.
Since $\eta_F=dF'(\eta)$ and $\eta_{\hat{F}\circ i}=d\hat F(\eta)$ belong
to $TN_0$, it is easily seen that $T_1$ is also a parallel vector bundle
isometry with ${\Delta}=\N(\alpha^{F}_{\Les_1^\perp})\cap
\N(\alpha^{\hat F\circ i}_{\hat{\Les}_1^\perp})$ and
$\alpha^{\hat{F}\circ i}_{\hat{\Les}_1}=T_1\circ \alpha^{F}_{{\Les}_1}$.
Moreover, since $\Les\subset \Les_1$, the restriction
$T=T_1|_\Les \colon\, \Les\to \hat{\Les}$
defines a parallel vector bundle isometry such that
$\beta^{\hat{F}\circ i}_{\hat{\Les}}=T\circ \beta^{F}_{{\Les}}$.

By \eqref{e:sffs3} we have
$$
P\circ \beta^{\hat{F}\circ i}= d\Psi(\va\beta^{\bar F}),
$$
where
$P\colon\,T_{\hat{F}\circ i}\Le^{n+q+2}\to d\Psi(T_{\bar F}\R^{n+q})$
denotes the orthogonal projection. This implies, in particular, that
$\Tes$ is well defined, and that $\Delta=
\N(\beta^{F}_{\Les^\perp})\cap \N(\beta^{\bar F}_{\bar\Les^\perp})$.
In addition,
$$
d\Psi\circ\Tes=P\circ T,
$$
which implies by \eqref{e:sffs3} that $\Tes$ is a vector bundle isometry.
Moreover,
$$
d\Psi\circ \Tes\circ \beta^F_{\Les}=P\circ T\circ \beta^F_{\Les}=P\circ
\beta^{\hat{F}\circ i}_{\hat\Les}=d\Psi(\va\beta^{\bar{F}}_{\bar\Les}),
$$
and therefore $ \Tes\circ\beta^{F}_{\Les}=\va\beta^{\bar{F}}_{\bar \Les}$
as desired.

Finally, we must prove that $\Tes$ is parallel with respect to the induced
connections. For $\xi\in \Les$ and $Y\in TN$, we have
\begin{eqnarray*}
d\Psi(\Tes (\nabla^\perp_Y \xi)_{{\Les}})&=&
P(T(\nabla^\perp_Y \xi)_{\bar\Les})=P(\nabla^\perp_Y T\xi)_{\hat\Les}=
P(\tilde\nabla_Y T\xi)_{\hat\Les}=(P\tilde\nabla_Y T\xi)_{P\hat\Les}\\
&=&(\tilde\nabla_Y PT\xi)_{P\hat\Les}=
(\tilde\nabla_Y d\Psi\Tes \xi)_{d\Psi\bar\Les}=
(d\Psi\tilde\nabla_Y \Tes \xi)_{d\Psi\bar\Les}=
d\Psi(\tilde\nabla_Y \Tes \xi)_{\bar\Les}\\
&=&d\Psi(\nabla^\perp_Y \Tes \xi)_{\bar\Les}\ ,
\end{eqnarray*}
and the claim follows.

\vskip .2cm
\noindent {\it $ii)$ The pair $\{f,\hat f\}$ is degenerate.}
\vskip .2cm

By \pref{deg1}, the pair $\{f',\hat f\}$ extends isometrically to
mutually $\Delta^{s_0}_0$-ruled Loren\-tz\-ian cones
$F_0'\colon\, N_0^{n+r_0}\to\Le^{n+p+2}$ and
$\hat F_0\colon\, N_0^{n+r_0}\to\Le^{n+q+2}$, with
\begin{equation}\label{e:est}
s_0\geq n-p-q+3\ell-4
\end{equation}
and $2\leq r_0\leq \ell$. Moreover, we have a parallel vector bundle
isometry $T_0:\Les_0^{\ell_0}\to\hat\Les_0^{\ell_0}$ that preserves
second fundamental forms, with $\ell=\ell_0+r_0$ and $\Delta_0^{s_0}=
\N(\a^{F_0'}_{\Les_0^\perp})\cap \N(\a^{\hat F_0}_{\hat\Les_0^\perp})$.

Since $f'$ is tangent to $F_0'$, $e_0$ is nowhere normal to
$F_0'$. Thus, $F_0'$ is transversal to the degenerate hyperplane
${\cal H}={\cal H}^{n+p+1}:=\{x\in\Le^{n+p+2}:\<x,e_0\>=1\}$, and
we locally define $N_1^{n+r_0-1}:=F_0'^{-1}(F_0'(N_0)\cap {\cal
H})\subset N_0^{n+r_0}$,  $\Delta_1^{s_0-1}:=\Delta^{s_0}_0\cap
TN_1$, and $F_1':=F'_0|_{N_1}$, $\hat F_1:=\hat F_0|_{N_1}$.

Now, $F_1'$ is transversal to $\Ve^{n+p+1}$, hence we may locally define
$N^{n+r} := F_1'^{-1}(F_1'(N_1)\cap \Ve^{n+p+1})\subset N_1^{n+r_0-1}$,
$\Delta^s:=\Delta_1^{s_0-1}\cap TN$, and $F':=F'_1|_N$,
$\hat F:=\hat F_1|_N$, with $s=s_0-2$, $r=r_0-2$.

Since $F'(N)\subset {\cal H}\cap \Ve^{n+p+1}=\Ee^{n+p}$, there exists
$F\colon\,N\to \R^{n+p}$ such that $F'=\Psi\circ F$. On the other
hand, using that $\<F_0',F_0'\>=\<\hat{F}_0,\hat{F}_0\>$, it follows that
$\hat F$ takes values in $\Ve^{n+q+1}$, and we may define
$\bar{F}= {\cal C}(\hat{F})\colon\,N\to \R^{n+q}$.
Then, as in the nondegenerate case, we obtain that $\{F,\bar{F}\}$ is a
conformal pair, $F\circ j=f$ and $\bar{F}\circ j=\bar{f}$, where $j$ is
the inclusion of $M$ into $N$, hence $\{F,\bar{F}\}$ is a conformal
extension of $\{f,\bar f\}$. Moreover, $F$ and $\bar{F}$ are mutually
$\Delta^s$-conformally ruled.

The estimate on $s$ also follows as in the nondegenerate case. From
\eqref{e:est} we have
$$
s=s_0-2\geq n-p-q+3(\ell_0+r),
$$
so it suffices to show that $\ell_0\geq \ell^c$. As before,
$\Les':=L^c_{\Delta}({F}')\subset {\Les}_0\oplus \spa\{dF_1'(\eta)\}$,
where $\eta(x)\in T_xN_1 \cap T_x^\perp N$ is the normal component of the
mean curvature vector at $x\in N$ of the leaf of $\Delta^s$ through $x$.
On the other hand,  $d\Psi(\Les)\subset \Les'\oplus \spa\{{F}'\}$, hence
\begin{equation}\label{e:lor2}
d\Psi(\Les)\subset \Les_0\oplus \spa\{dF_1'(\eta)\}\oplus \spa\{F'\}.
\end{equation}
Arguing as before, we obtain that the subspace $W$ on the right-hand-side
of \eqref{e:lor2} is Lorentzian. Therefore, $d\Psi(\Les)$ is a Riemannian
subspace of $W$ that is orthogonal to the null vector ${ F}'\in W$, hence
it has codimension at least two in $W$.

We claim that \eqref{e:iso}  defines also in this case a parallel
vector bundle isometry $\Tes\colon\, \Les\to\bar\Les$
such that $\Tes\circ\beta^{F}_{\Les}=\va\beta^{\bar{F}}_{\bar \Les}$.

Choose smooth unit vector fields $\xi'$ and $\hat{\xi}$ spanning
$T_{F_0'}N_0\cap T^\perp_{F_1'}N_1$ and
$T_{\hat{F}_0}N_0\cap T^\perp_{\hat{F}_1}N_1$, respectively, and set
${\Les}_1:={\Les}_0\oplus \spa\{\xi'\}\subset T_{F_1}^\perp N_1$ and
$\hat{\Les}_1:=
\hat{\Les}_0\oplus \spa\{\hat{\xi}\}\subset T_{\hat{F}_1}^\perp N_1$.
Extend the parallel vector bundle isometry
$T_0\colon\,\Les_0\to \hat{\Les}_0$ to
$T_1\colon\,\Les_1\to \hat{\Les}_1$ by setting
$T_1|_{{\Les}_0}=T_0$ and $T_1(\xi')=\hat{\xi}$.

Now set $\Les_2:=\Les_1\oplus\spa\{\eta_{F'}\}\subset T_{F'}^\perp N$ and
$\hat{\Les}_2:=\hat\Les_1\oplus \spa\{\eta_{\hat F}\}\subset
T_{\hat{F}}^\perp N$ and extend $T_1$ to
$T_2\colon\,\Les_2\to \hat{\Les}_2$ by setting
$T_2|_{{\Les}_1}=T_1$ and $T_2(\eta_{F'})=\eta_{\hat F}$. Since
$\eta_{F'}=dF_1'(\eta)$ and $\eta_{\hat F}=d\hat F_1(\eta)$ belong to
$TN_1$, it is easily seen that $T_2$ is also a parallel vector bundle
isometry with
$\Delta^s=\N(\alpha^{F'}_{\Les_2^\perp})\cap
\N(\alpha^{\hat F}_{\hat\Les_2^\perp})$
and
$\alpha^{\hat{F}}_{\hat{\Les}_2}=T_2\circ \alpha^{F'}_{{\Les}_2}$.
Moreover, since $\Les\subset \Les_2$, the restriction
$T=T_2|_\Les \colon\, \Les\to \hat{\Les}$
defines a parallel vector bundle isometry
such that $\beta^{\hat{F}}_{\hat{\Les}}=T\circ \beta^{F'}_{{\Les'}}$.

By \eqref{e:sffs3} we have
$$
P\circ \beta^{\hat{F}}= d\Psi(\va\beta^{\bar F})\,\,\,\mbox{and}\,\,\,
\beta^{F'}=d\Psi\circ \beta^F.
$$
This implies, in particular, that $\Tes$ is well defined and that $\Delta=
\N(\beta^{F}_{\Les^\perp})\cap \N(\beta^{\bar F}_{\bar\Les^\perp})$.
Moreover,
$$
d\Psi\circ\Tes=P\circ T\circ d\Psi,
$$
which implies by \eqref{e:sffs3} that $\Tes$ is a vector bundle isometry.
Thus,
$$
d\Psi\circ \Tes\circ \beta^F_{\Les}=
P\circ T\circ d\Psi\circ \beta^F_{\Les}=
P\circ \beta^{\hat{F}}_{\hat\Les}=d\Psi(\va\beta^{\bar{F}}_{\bar\Les}),
$$
and the claim is proved.

The proof that $\Tes$ is parallel with respect to the induced
connections is also analogous to that of the nondegenerate case.
\qed

\section{The proofs of \tref{int2} and \cref{cor1}}

\noindent{\em Proof of \tref{int2}:\/} If the pair $\{f,\hat{f}\}$ in the
proof of \tref{int3} is nondegenerate, take the maximal isometric
$\Delta^{d+r}$-ruled extensions $F'\colon\,N^{n+r}\to\R^{n+p}$ and
$\hat F\colon\,N^{n+r}\to\Le^{n+q+2}$, $r\geq 1$, of $f$ and $\hat f$, respectively,
given by Theorem \ref{t3}.

If $\{f,\hat{f}\}$ is a degenerate pair, let
$F_1'\colon\,N_1^{n+r}\to {\cal H}\subset \Le^{n+p+2}$ and
$\hat F_1\colon\,N_1^{n+r}\to \Le^{n+q+2}$ be as in the proof of
\tref{int3} for the degenerate case. Notice that $r\geq 1$, since
$N_1^{n+r}$ has dimension one less
than that of the extension $F_0'\colon\,N_0^{n+r_0}\to \Le^{n+q+2}$ of
$f'$ as in the same proof, and $r_0\geq 2$.

We claim that $F':=\pi\circ F_1'$ is an immersion, where
$\pi:\Le^{n+p+2}=\Le^2\times\R^{n+p}\to\R^{n+p}$ is the projection
onto the second factor. Otherwise, $N_1^{n+r}$ is degenerate or,
equivalently, $e_0\in T_{F_1'}N_1\subset
T_{F_0'}N_0=T_{f'}M\oplus\Lambda$ (see \eqref{e:par}). That is,
there is $X_0\in TM$ such that \mbox{$X_0+e_0\in \Delta_0 =
\N(\phi)$} as in \eqref{bili}. By \eqref{e:ssf} and
$\Ies(e_0)=\xi_0$, this implies that \mbox{$\Ima d\xi_0 = \Ima
\hat\nabla^\perp \xi_0 \subseteq L$}. On the other hand, $e_0$ is
also normal to $f'$, so $e_0\in L \subseteq S_0$. It follows that
$\xi_0$ is constant, since $d\xi_0=\hat\nabla^\perp
\xi_0=\nabla^\perp e_0=0$.

We conclude that $\hat f=\bar\Psi\circ f^*$ for some isometric
immersion $f^*\colon\,M^n\to \R^{n+q}$, where $\bar
\Psi\colon\,\R^{n+q}\to \Ve^{n+q+1}\subset \Le^{n+q+2}$ is the
isometric embedding defined as in \eqref{e:psi} with respect to a
pseudo-orthonormal basis of $\Le^{n+q+2}$ containing $\xi_0$.
Since $\bar \Psi=T\circ \Psi$ for some orthogonal linear
transformation $T$ of $\Le^{n+q+2}$, it follows that ${\cal
I}(\bar{f})=\hat f=\bar\Psi\circ f^*=T\circ \Psi\circ f^*= T\circ
{\cal I}({f}^*)$, that is, ${\cal I}(\bar{f})$ and ${\cal
I}({f}^*)$ are isometrically congruent in $\Le^{n+q+2}$. Therefore
$\bar f$ and $f^*$ are conformally congruent in $\R^{n+q}$. This
contradicts our hypothesis and proves our claim.

Thus, setting $N=N_1$, also in the degenerate case we obtain maximal
isometric \mbox{$\Delta^{d+r}$-ruled} extensions
$F'\colon\,N^{n+r}\to\R^{n+p}$ and
$\hat F=\hat F_1\colon\,N^{n+r}\to\Le^{n+q+2}$, $r\geq 1$, of $f$ and
$\hat f$, respectively.

As in the proof of \tref{int3}, since $\hat{F}$ is transversal to the
light cone, by restricting to an open subset, if necessary, we may assume
that $\hat{F}$ is an embedding, so that
\mbox{$\bar N=\hat{F}^{-1}(\hat{F}(N)\cap \Ve^{n+q+1})\supset M^n$} is an
$(n+r-1)$-dimensional manifold. As before, setting $F=F'\circ i$ and
$\bar{F}= {\cal C}(\hat{F}\circ i)\colon\,\bar{N}\to \R^{n+q}$, where
$i\colon\, \bar N\to N$ is the inclusion map, we have that
$\{F,\bar{F}\}$ is a conformal pair, $F\circ j=f$ and
$\bar F\circ j=\bar f$, where $j$ is the inclusion of $M$ into $\bar N$.
Thus, $\{F,\bar{F}\}$ is a conformal extension of $\{f,\bar f\}$.

Since $\{f,\bar f\}$ is a genuine conformal pair, we must have $r=1$,
hence $\bar N=M$, $F\circ i=f$ and ${\cal C}(\hat{F}\circ i)=\bar{f}$.
A similar argument shows that any isometric extension of the pair
$\{F',\hat F\}$ would give a conformal extension of the pair
$\{f,\bar f\}$, hence $\{F',\hat F\}$ must be a genuine isometric
pair.\qed\vspace{1.5ex}

\noindent {\em Proof of \cref{cor1}:\/}
Given a conformal immersion $\bar f\colon\,M^n\to \R^{n+q}$, \tref{int3}
applies and yields, locally on an open dense subset of $M^n$, a (possibly
trivial) conformal extension $\{F,\bar F\}$ of $\{f,\bar f\}$ to a
mutually $\Delta^s$-conformally ruled pair of immersions
$F\colon N^{n+r}\!\to\R^{n+p}$ and $\bar F\colon N^{n+r}\!\to\R^{n+q}$,
$0\leq r\leq p$, with
$s\geq n-p-q+3(\ell^c+r)$. It suffices to prove that $r=p$.

First we show that $L^c_{\Delta}(f)^\perp=\{0\}$. Otherwise, if
$s':=\rank L^c_{\Delta}(f)^\perp>0$, we would have, since
$\Delta\subset \N(\beta^f_{{L^c_{\Delta}}^\perp})$, that
$$
\nu^c_{s'}\geq s\geq n-p-q+3(p-s') = n+p-q-2s'+(p-s'),
$$
contradicting our assumption on $\nu^c_{s'}$. Therefore, we have
$$
s\geq n+2p-q.
$$

Now assume that $r<p$. Since $F$ is $\Delta$-conformally ruled, we have
that $\alpha_F=\<\ ,\ \>\eta_F$ on $\Delta\times\Delta$
for some normal vector field $\eta_F$. In particular, for any unit normal
vector field $\xi\in T_F^\perp  N$ we obtain that
$\<(A^f_\xi -\<\xi,\eta_F\>I) D,D\>=0$, where $D= \Delta\cap TM$.
Since $\rank D =s-r\geq n-(q-p)+1$, then
$\nu^c_1 \geq n-2(q-p)+2$, and this is a contradiction with our
assumption on $\nu^c_1$.
\qed

\vspace{.5in} {\renewcommand{\baselinestretch}{1}

\hspace*{-20ex}\begin{tabbing} \indent     \= IMPA -- Estrada Dona
Castorina, 110
\indent\indent  \=  Universidade Federal de S\~{a}o Carlos \\
\> 22460-320 -- Rio de Janeiro -- Brazil  \>
13565-905 -- S\~{a}o Carlos -- Brazil \\
\> E-mail: luis@impa.br \> E-mail: tojeiro@dm.ufscar.br
\end{tabbing}}
\end{document}